\documentstyle[12pt,twoside]{article}
\newtheorem{theorem}{Theorem}

\newtheorem{prop}[theorem]{Proposition}
\newtheorem{lemma}[theorem]{Lemma}

\newtheorem{definition}[theorem]{Definition}
\newtheorem{conjecture}[theorem]{Conjecture}

\begin{document}
%\begin{document}
\overfullrule=0pt
\baselineskip=24pt
% if i uncomment the preceeding line, the text will be double-spaced.
\font\tfont= cmbx10 scaled \magstep3
\font\sfont= cmbx10 scaled \magstep2
\font\afont= cmcsc10 scaled \magstep2
\title{\tfont On periodic sequences for algebraic numbers}
\bigskip
\author{ Thomas Garrity \\
\\ Department of Mathematics\\ Williams College\\Williamstown, MA  
01267\\ 
email:tgarrity@williams.edu}
\date{}
\maketitle
\begin{abstract}
% Abstract text begins here

 For each positive integer $n \geq 2$, a 
 new approach to expressing real numbers as sequences of nonnegative 
integers is given. The $n=2$ case is equivalent to 
the standard continued fraction 
algorithm.  For $n=3$, it reduces to a new iteration of the 
triangle.  Cubic irrationals that are roots of $x^{3} + kx^{2} + x - 
1$ are shown to be precisely those numbers with purely periodic 
expansions of period length one.  For  general positive integers $n$, 
it reduces to a new iteration of an $n$ 
dimensional simplex. 
\end{abstract}

\section{Introduction}

In 1848 Hermite \cite{Herm} posed to Jacobi the problem of 
generalizing continued fractions so that periodic expansions of a 
number reflect its algebraic properties.  We state this as:

\noindent{\bf The Hermite Problem}: {\it Find methods for writing 
	numbers that reflect special algebraic properties.}

\noindent	In attempting to answer this question, Jacobi developed a special case 
	of what is now called the Jacobi-Perron algorithm. 
	Bernstein \cite{Bernstein} wrote a good survey of this algorithm; 
	Schweiger 
   \cite{Schweiger} covered its ergodic properties. Brentjes' book \cite{Brentjes} is 
  a good source for its many variations. Using quite different methods,  
  Minkowski  \cite{Minkowski 1899}  developed a quite different 
  approach to the 
	 Hermite's problem. For another attempt, see the work of Ferguson 
	 and Forcade
	 \cite{Ferguson-Forcade}.

In this paper we give another approach, 
which will also be a  generalization of continued fractions.  
  To each $n$-tuple of real numbers
 $(\alpha_{1}, \ldots ,\alpha_{n})$, 
 with $1 \geq \alpha_{1} \geq  \ldots \geq \alpha_{n}$, we
  will associate  a sequence 
 of nonnegative integers.  For reasons that will become apparent later, we 
 will call this sequence the {\it triangle sequence} (or {\it simplex 
 sequence}) for the $n$-tuple. The hope is that the periodicity of 
 this sequence will provide insight into whether or not 
  the $\alpha_{k}$ are algebraic of degree at most $n$.  We will show 
  that this is the case for when $n=3$.
% The pairs of particular interest will be 
% 	 $(\alpha, \alpha^{2})$.  If the triangle sequence terminates, we will see 
% 	 that this means that there are integers $p$, $q$ and $r$, not all 
% 	 zero, such that 
% 	 $$p + q \alpha + r \beta = 0.$$
% 	 This will be seen to be the analogue of the fact that the sequence 
% 	 for continued 
% 	 fractions terminates precisely when the initial real is rational.  
% 	 We will also see that if the triangle sequence is eventually periodic, then 
%  $\alpha$ and $\beta$ are at worst  cubic irrationals and $\alpha$
%  and $\beta$ are algebraically related.  Further, if the 
%  triangle sequence is periodic with period one, say with sequence 
%  $(a,a,a, \ldots)$, then $\beta = \alpha^{2}$ and more importantly
%   $\alpha$ is a root of the cubic 
%  $$x^{3} + a x^{2} + x - 1 = 0,$$
%  in direct analogue to the fact that numbers with periodic continued 
%  fraction expansion with period length one are roots of $x^{2} - a x - 
%  1 =0$.
%  
%  One of the main justifications for this approach is that it involves 
%  a quite interesting iteration of the triangle
%  $$ \bigtriangleup = \{ (x,y): 1 \geq x \geq y \geq 0 \} .$$
%  This triangle will be the analogue of the unit interval for 
%  continued fractions.

In the next section we quickly review some well-known facts about
 continued fractions.  We then concentrate on the cubic case, for 
 ease of exposition.  The proofs go over easily to the general case, 
 which we will discuss in  section nine.  In section three we define, 
 given a pair $(\alpha, \beta) \in \{ (x,y): 1 \geq x \geq y \geq 0 \}$,
  the triangle 
 iteration and the triangle sequence.  Section four will recast the 
 triangle sequence via matrices.  This will allow us to interpret the 
 triangle sequence as a method  for producing integer lattice points 
 that approach the 
 plane $x + \alpha y + \beta z=0$.  Section five will show that 
 nonterminating triangle sequences  uniquely determine the pair $(\alpha, 
 \beta)$.  Section six discusses how  every possible triangle 
 sequence corresponds to
  a pair $(\alpha, \beta) \in \{ (x,y): 1 \geq x \geq y \geq 0 \}.$
  Section seven turns to classifying those pairs with purely periodic 
  sequences.  Section eight concerns itself with periodicity in 
  general.  Section nine deals with the general case $n$. 
  
  At http://www.williams.edu/Mathematics/tgarrity/triangle.html, 
  there is a web page that gives many examples of triangle sequences 
  and provides software packages running on Mathematica for making 
  actual computations.

I would like to thank Edward Burger for many helpful discussions.  
Tegan Cheslack-Postava, Alexander 
Diesl, Matthew Lepinski and  Adam Schuyler  have provided critical aid.
 I 
would also like to thank K. M. Briggs for useful comments on an 
earlier version of this paper.
 \section{Continued Fractions}
 Given a real number $\alpha$, recall that its continued fraction 
 expansion is:
 $$	 \alpha  = a_{0} + {1\over {a_{1}  + {1\over a_{2}  + ...}}},$$
where  
$		a _{0} = \{ \alpha \} = \mbox{greatest integer part of}\; \alpha$,
$$a_{1} = \{ {1\over \alpha  - a_{0}} \}\; \mbox{ and}\;
		b_{1} =  {1\over \alpha  - a_{0}} - a_{1}.$$
Inductively, define 
$$		a_{k} = \{{1\over b_{k-1}}\} \;
\mbox{and}\;
b_{k} = {1\over b_{k-1}} - 	a_{k}.$$
A number's continued fraction expansion can be captured by 
examining iterations of the Gauss map $G:I \rightarrow I$, with $I$ denoting 
the unit interval $(0,1]$, defined by 
$$G(x) = {1\over x} - \{{1\over x}\}.$$
If we partition the unit interval into a disjoint union of 
subintervals:
$$I_{k} = \{ x \in I: \frac{1}{k+1} < x \leq \frac{1}{k} \},$$
then the nonnegative integers $a_{k}$ in the continued fraction 
expansion of $\alpha$ can be interpreted as keeping track of which 
subinterval the number $\alpha$ maps into under the kth iterate of $G$. 
Namely,  $G^{k}(\alpha)\in I_{a_{k}}$.

\section{The Triangle Iteration}
In this section we define an iteration $T$ on the triangle
$$ \bigtriangleup = \{ (x,y): 1 \geq x \geq y >  0 \}.$$
Partition this triangle into disjoint triangles
$$\bigtriangleup_{k} = \{ (x,y) \in \bigtriangleup: 1 - x - ky 
\geq 0 > 1 - x - (k+1)y \},$$
where $k$ can be any nonnegative integer. Note that its vertices are 
$(1,0)$, $(\frac{1}{k+1},\frac{1}{k+1})$ and $(\frac{1}{k+2},\frac{1}{k+2})$.

 Define the {\it triangle map}
 $T:\bigtriangleup \rightarrow \bigtriangleup \cup \{ (x,0): 0 \leq 
 x\leq 1 \}$ by 
setting
$$T(\alpha,\beta) = (\frac{\beta}{\alpha}, \frac{1 - \alpha - k 
\beta}{\alpha}),$$
if the pair $(\alpha,\beta) \in \bigtriangleup_{k}$. Frequently we 
will abuse notation by denoting $\bigtriangleup \cup \{ (x,0): 0 \leq 
 x\leq 1 \}$ by  $\bigtriangleup$.

We want to associate a sequence of nonnegative integers to the 
iterates of the map $T$.  Basically, if $T^{k}(\alpha,\beta) \in 
 \bigtriangleup_{a_{k}}$, we will associate to $(\alpha, \beta)$ 
 the sequence $(a_{1}, \ldots)$.

Recursively define a sequence of decreasing positive reals and a sequence of 
nonnegative integers as follows:
Set $d_{-2}= 1, d_{-1} = \alpha, d_{0}=\beta.$ Assuming that we have 
$d_{k-3} > d_{k-2}> d_{k-1}>0$, define $a_{k}$ to be a nonnegative 
integer such that
$$d_{k-3} -d_{k-2} - a_{k} d_{k-1} \geq 0$$
but
$$d_{k-3} -d_{k-2} - (a_{k} + 1) d_{k-1} < 0.$$
Then set
$$d_{k} = d_{k-3} -d_{k-2} - a_{k} d_{k-1}.$$
If at any stage $d_{k}=0$, stop.  

\begin{definition}
The triangle sequence of the pair $(\alpha, \beta)$ is the sequence
$(a_{1}, \ldots)$.
\end{definition}
We will say that the triangle sequence {\it terminates} if at any 
stage $d_{k}=0$.  In these cases, the triangle sequence will be 
finite.

Note that 
$$T(\frac{d_{k-1}}{d_{k-2}}, \frac{d_{k}}{d_{k-2}})
= (\frac{d_{k}}{d_{k-1}}, \frac{d_{k+1}}{d_{k-1}}).$$

Also note that by comparing this to the first part of chapter seven 
in \cite{Stark}, we see that this is indeed a generalization of 
continued fractions.

\section{The Triangle Iteration via Matrices and Integer Lattice 
Points}

Let $(a_{1}, \ldots)$ be a triangle sequence associated to the 
pair $(\alpha, \beta)$.  Set
$$P_{k} = \pmatrix{0&0&1\cr 1&0&-1 \cr 0&1&-a_{k}\cr}.$$
Note that $\det P_{k} =  1$.
Set $M_{k} = P_{1}\cdot P_{2} \cdots P_{k}$. This allows us to 
translate the fact that 
$$T(\frac{d_{k-1}}{d_{k-2}}, \frac{d_{k}}{d_{k-2}})
= (\frac{d_{k}}{d_{k-1}}, \frac{d_{k+1}}{d_{k-1}})$$
 into the language 
of matrices via the following proposition (whose proof is 
straightforward):
\begin{prop}
Given the pair $(\alpha, \beta)$, we have
$$(d_{k-2}, d_{k-1}, d_{k}) = (1, \alpha, \beta)M_{k}.$$
\end{prop}
 Write 
$$M_{k} = \pmatrix{p_{k-2}&p_{k-1}&p_{k}\cr q_{k-2}&q_{k-1}&q_{k}
 \cr r_{k-2}&r_{k-1}&r_{k}\cr}.$$
 Then a calculation leads to:
 \begin{prop}
 For all $k$, we have
 $$p_{k} = p_{k-3} - p_{k-2} - a_{k}p_{k-1},$$
 $$q_{k} = q_{k-3} - q_{k-2} - a_{k}q_{k-1},$$
and
 $$r_{k} = r_{k-3} - r_{k-2} - a_{k}r_{k-1}.$$
\end{prop}

Set
$$C_{k} = \pmatrix{p_{k}\cr q_{k}
 \cr r_{k}\cr}.$$
 Note that $C_{k}$ can be viewed as a vector in the integer lattice.  
Then the numbers $d_{k}$ are seen to be a measure of the distance from the 
 plane $x + \alpha y + \beta z = 0$ to the lattice point $C_{k}$, 
 since $d_{k} = (1, \alpha, \beta)C_{k}.$  Observing that 
 $$C_{k}= C_{k-3} - C_{k-2} - a_{k}C_{k-1},$$
 we see thus that the triangle sequence encodes information of how to get a 
 sequence of lattice points to approach the plane $x + \alpha y + \beta 
 z = 0$, in direct analogue to continued fractions \cite{Stark}.  
Unlike the continued fraction case, though, these lattice points need 
not be the best such approximations.

\section{Arbitrary triangle sequences}

\begin{theorem}
Let $(k_{1}, k_{2}, 
\ldots )$ be any infinite sequence of nonnegative integers with 
infinitely many of the $k_{i}$ not zero.  Then there is a  pair 
$(\alpha, \beta)$ in $\bigtriangleup$ that has this sequence as its 
triangle sequence.
\end{theorem}

\noindent {\bf Proof:}
Suppose that we have an infinite triangle sequence $(k_{1}, k_{2}, 
\ldots )$.  
By a straighforward calculation,we see that a line with 
equations $y=mx+b$ will map to the line
$$(1-k b)u - (1- k b)m = bv + bkm + b,$$
where $T(x,y) = (u,v).$

The map $T$, restricted to the triangle $\bigtriangleup_{k}$,
 will send  the vertices of 
$\bigtriangleup_{k}$ to the vertices of $\bigtriangleup$, with 
$T(1,0) = (0,0)$, $T(\frac{1}{k+1}, \frac{1}{k+1})= (1,0)$ 
and  $T(\frac{1}{k+2}, \frac{1}{k+2})= (1,1)$.
 Restricted to $\bigtriangleup_{k}$, the map $T$ is thus one-to-one and 
 onto $\bigtriangleup$.

 But this gives us our theorem, as each $\bigtriangleup_{k}$ 
 can be split into its own (smaller) triangles, one for each 
 nonnegative integer, and hence each of these smaller triangles can be 
 split into even smaller triangles, etc. Hence to each nonterminating 
 triangle sequence there corresponds a pair $(\alpha, \beta)$.  
 \noindent {\bf QED}

 \section{Recovering points from the triangle sequence}

 The question of when a triangle sequence determines a 
 unique pair $(\alpha, \beta)$ is subtle.  If the sequence 
 terminates, then  the pair $(\alpha, \beta)$ is not unique.
 Even if the triangle sequence does not terminate, we do not 
 necessarily have uniqueness, as discussed in \cite{CDGLS}.
 But we do have 
 \begin{theorem}
 If an integer $k$ occurs infinitely often in a
    nonterminating  sequence $(k_{1}, k_{2}, 
\ldots )$ of nonnegative integers, then
 there is a unique pair $(\alpha, \beta)$ in $\bigtriangleup$
that has this sequence as its triangle sequence.
\end{theorem}
The proof is contained in \cite{CDGLS} and is not easy.
 
 If the triangle sequence uniquely determines a pair $(\alpha, \beta)$,
 then we can recover $(\alpha, \beta)$ as follows. 
  By construction, the numbers $d_{k}$  approach 
 zero.  Consider the plane
 $$x + \alpha y + \beta z = 0,$$
 whose normal vector is $(1, \alpha, \beta)$.  As seen in the last section,
  the columns of the 
 matrices $M_{k}$ can be interpreted as vectors that are approaching 
 this plane. This will allow us to prove: 
 \begin{theorem}
 If a triangle sequence uniquely determines  the pair $(\alpha, \beta)$, then 
 $$\alpha = \lim_{k\rightarrow \infty} \frac{p_{k}r_{k-1} - 
 p_{k-1}r_{k}}{q_{k-1}r_{k}-q_{k}r_{k-1}}$$
 and
 $$\beta = \lim_{k\rightarrow \infty}\frac{p_{k-1}q_{k} - 
 p_{k}q_{k-1}}{q_{k-1}r_{k}-q_{k}r_{k-1}}.$$
 \end{theorem}
 The proof is also in \cite{CDGLS}.  The quick, but incorrect, argument is 
 that 
  the vectors $(p_{k-1},q_{k-1},r_{k-1})$ and 
 $(p_{k},q_{k},r_{k})$ are columns in the matrix $M_{k}$, each
 of which  approaches being in the plane  $x + \alpha y + \beta z = 
 0$.  Thus the limit as $k$ approaches infinity of the 
 cross product of these two vectors must point in the 
 normal direction  $(1, \alpha, \beta)$.  But this is the above limits.

\section{Purely periodic triangle sequences of period length one}
\begin{theorem}
Let $0 < \beta \leq \alpha < 1$ be a pair of numbers whose triangle 
sequence is $(k,k,k,\ldots)$.  Then $ \beta = \alpha^{2}$ and $\alpha$
is a root of the cubic equation
$$x^{3} + k x^{2} + x -1 = 0.$$
Further if $\alpha$ is the real root of this cubic that is between
zero and one, then 
$(\alpha, \alpha^{2})$ has purely periodic triangle sequence $(k,k,k, 
\ldots)$. 
\end{theorem}
% NEED TO HAVE PROOF OF UNIQUENESS OF TRIANGLE SEQUENCE BEFORE THE 
% NEXT STATEMENT CAN BE TRUE  
\noindent {\bf Proof:} We need $T(\alpha,\beta) = (\alpha,\beta)$. Since $T(\alpha,\beta)= 
 (\frac{\beta}{\alpha}, \frac{1-\alpha - k \beta}{\alpha})$, we need 
 $$\alpha = \frac{\beta}{\alpha}$$ 
 and
 $$\beta = \frac{1-\alpha - k \beta}{\alpha}.$$
 From the first equation we get $\beta = \alpha^{2}$.  Plugging in 
 $\alpha^{2}$ for $\beta$ in the second equation and clearing 
 denominators leads to 
 $$\alpha^{3} + k\alpha^{2} + \alpha -1 = 0$$
 and the first part of the theorem.
 
 Now for the converse.  Since the polynomial $x^{3} + k x^{2} + x -1 $ 
 is $-1$ at $x=0$ and is positive at $x=1$, there is root $\alpha$ 
 between zero and one.  We must show that $(\alpha, \alpha^{2})$ is 
 in $\bigtriangleup_{k}$ and that $T(\alpha, \alpha^{2}) = (\alpha, \alpha^{2})$.
We know that 
$$\alpha^{3} = 1 - \alpha - k \alpha^{2}.$$
Since $\alpha^{3}>0$, we have $1 - \alpha - k \alpha^{2} > 0$.
Now
$$ 1 - \alpha - (k+1) \alpha^{2} = \alpha^{3} - \alpha^{2} < 0,$$
which shows that $(\alpha, \alpha^{2}) \in \bigtriangleup_{k}$.

Finally,
$$T(\alpha, \alpha^{2}) = (\frac{\alpha^{2}}{\alpha}, 
\frac{1-\alpha - k \alpha^{2}}{\alpha})$$
$$= (\alpha, \alpha^{2}).$$

\noindent {\bf QED}

Similar formulas for purely periodic sequences with period 
length two, three, etc., can be computed, but they quickly 
become computationally 
messy.

 \section{Terminating and Periodic Triangle Sequences}
 
 We first want to show that if $(\alpha, \beta)$ is a pair of rational 
 numbers, then the corresponding triangle sequence must terminate, 
 meaning that eventually all of the $k_{n}$ will be zero.
 \begin{theorem}
 Let $(\alpha, \beta)$ be a pair of rational numbers in 
 $\bigtriangleup$.  Then the corresponding triangle sequence terminates.
 \end{theorem}
 
 \noindent {\bf Proof:}
 In constructing the triangle sequence, we are 
 just concerned with the ratios of the triple $(1, \alpha, \beta)$.  
 By clearing denominators, we can replace this triple by a triple of 
 positive integers $(p,q,r)$, with $p \geq q \geq r$. Then we have 
 $d_{-2} = p, d_{-1}=q, d_{0}=r$.  Then the sequence of $d_{k}$ will 
 be a sequence of positive decreasing integers.  Thus for some $k$ we 
 must have $d_{k}=0$, forcing the triangle sequence to terminate. 
 
\noindent {\bf QED}

 Now to see what happens when the triangle sequence is eventually 
 periodic.   
 \begin{theorem}
 Let $(\alpha, \beta)$ be a pair of real numbers in 
 $\bigtriangleup$ whose triangle sequence is eventually periodic.  
 Then $\alpha$ and  $\beta$ have degree at most three,  with 
 $\alpha \in {\bf Q}[\beta]$ or $\beta \in {\bf Q}[\alpha]$.
 \end{theorem}
 
 \noindent {\bf Proof:} If both $\alpha$ and $\beta$ are 
 rational, then by the above theorem the triangle sequence 
 terminates.  Thus we assume that not both $\alpha$ and and $\beta$ 
 are rational.  Since the triangle sequence is periodic, there will be 
 an integer appearing infinitely often in this sequence, which means 
 that the sequence will uniquely determine a pair $(\alpha, \beta)$.
 
 If the triangle sequence is periodic, there must be 
 an $n$ and $m$ so that 
 $$(\frac{d_{n-2}}{d_{n}}, \frac{d_{n-1}}{d_{n}})=
 (\frac{d_{m-2}}{d_{m}},\frac{d_{m-1}}{d_{m}}).$$
 Thus there exists a number $\lambda$ with 
 $$(d_{n-2},d_{n-1},d_{n})= \lambda (d_{m-2},d_{m-1},d_{m}).$$
 Using matrices we have:
 $$(1, \alpha, \beta)M_{n} = \lambda (1, \alpha, \beta)M_{m}$$
 and thus
 $$(1, \alpha, \beta)M_{n} M_{m}^{-1} = \lambda (1, \alpha, \beta).$$
  Since  $M_{n}$ and $M_{m}$ have integer
   coefficients, the matrix $M_{n} M_{m}^{-1}$ will have rational 
   coefficients.  Since the $d_{k}$ are decreasing, we must have 
   $|\lambda| \neq 1$. Since both $M_{n}$ and $M_{m}$ have 
   determinant  one, we have that $M_{n} M_{m}^{-1}$ 
   cannot be a multiple of the identity matrix.
   
   Set 
   $$M_{n} M_{m}^{-1}= \pmatrix{q_{11}&q_{12}&q_{13}\cr 
   q_{21}&q_{22}&q_{23}
 \cr q_{31}&q_{32}&q_{33}\cr}.$$
Then
$$q_{11} + q_{21}\alpha + q_{31} \beta = \lambda$$
$$q_{12} + q_{22}\alpha + q_{32} \beta = \lambda \alpha$$
and
$$q_{13} + q_{23}\alpha + q_{33} \beta = \lambda \beta.$$
We can eliminate the unknown $\lambda$ from the first and second equations 
and then from the first and third equations, leaving two equations 
with unknowns $\alpha$ and $\beta$.  Using these two equations we can 
eliminate one of the remaining variables, leaving the last as the 
solution to  polynomial with rational coefficients.  If this 
polynomial is the zero polynomial, then it can be seen that this will 
force $M_{n} M_{m}^{-1}$ to be a multiple of the identity, which is 
not possible.  Finally, it can be checked 
that this polynomial is a cubic.

\noindent {\bf QED}

\section{The higher degree case}

Almost all of this goes over in higher dimensions.  
We just replace our triangle by a dimension 
$n$ simplex.  The notation, though, is more cumbersome.

 Set
$$ \bigtriangleup = \{ (x_{1},\ldots ,x_{n}): 1 \geq x_{1} \geq \ldots
 \geq x_{n} > 0 \}.$$
 As we did before, we will frequently also call 
 $ \bigtriangleup = \{ (x_{1},\ldots ,x_{n}): 1 \geq x_{1} \geq \ldots
 \geq x_{n} \geq 0 \}.$
Set
$$\bigtriangleup_{k} = \{ (x_{1},\ldots ,x_{n})
 \in \bigtriangleup: 1 - x_{1} - \ldots  - x_{n-1} -  kx_{n} 
\geq 0 > 1 - x_{1} - \ldots  - x_{n-1} -  (k+1)x_{n}\},$$
where $k$ can be any nonnegative integer.  These are the direct 
analogue of the triangles $\bigtriangleup_{k}$ in the first part of 
this paper.  Unlike the earlier case, these $\bigtriangleup_{k}$, 
while disjoint, do not partition the simplex $\bigtriangleup$.  To 
partition $\bigtriangleup$, we need more simplices.
Set
$$\bigtriangleup ' = \{ (x_{1},\ldots ,x_{n})
 \in \bigtriangleup: 0 > 1 - x_{1} - \ldots  - x_{n} 
\}.$$
Then set
$$\bigtriangleup_{ij} = \{ (x_{1},\ldots ,x_{n})
 \in \bigtriangleup ': x_{j} \geq 1 - x_{1} - \ldots  - x_{i} 
\geq x_{j+1} \},$$
where $ 1 \leq i \leq n-2$ and $i < j \leq n$.  Also, we use the 
convention that $x_{n+1}$ is identically zero.
\begin{lemma}
The $\bigtriangleup_{k}$ and $\bigtriangleup_{ij}$ form a simplicial 
decomposition of the simplex $\bigtriangleup$.
\end{lemma}

\noindent {\bf Proof:}  It can be directly checked that the 
$\bigtriangleup_{k}$ and $\bigtriangleup_{ij}$ do form a disjoint 
partition of $\bigtriangleup$.  We need to show
 that the $\bigtriangleup_{k}$ and $\bigtriangleup_{ij}$
are simplices. Thus we want to show that each of these polygons have 
exactly $n+1$ vertices.  Label the $n+1$ 
vertices of the simplex $\bigtriangleup$ by $v_{0} = (0,\ldots, 0),
v_{1}=(1,0,\ldots, 0), v_{2}=(1,1,0,\ldots, 0), \ldots, 
v_{n}=(1,\ldots,1)$.  We label each of the $\frac{n(n+1)}{2}$ edges
 of the simplex by $v_{i}v_{j}$ if the endpoints of the edge are the 
 vertices $v_{i}$ and $v_{j}$. Consider 
the set  $\bigtriangleup_{k}$.   The hyperplanes 
$$x_{1}+\ldots+ x_{n-1} + k x_{n} = 1$$
and 
$$x_{1}+\ldots+ x_{n-1} + (k+1) x_{n} = 1$$
form two of the faces.  These hyperplanes intersect  each edge
$v_{0}v_{l}$, with $l < n$, in the same point $(\frac{1}{l}, \ldots, 
\frac{1}{l}, 0, \ldots,0)$, where this $n$-tuple has its first $l$ 
terms $\frac{1}{l}$ and the rest zero. The hyperplanes intersect the 
edge $v_{0}v_{n}$ in two distinct points: $x_{1}+\ldots+ x_{n-1} + k x_{n} = 1$
intersects at the point $(\frac{1}{n+k-1},\ldots,\frac{1}{n+k-1})$ 
while $x_{1}+\ldots+ x_{n-1} + (k+1) x_{n} = 1$ intersects in the point
$(\frac{1}{n+k},\ldots,\frac{1}{n+k})$. Since both hyperplanes contain the 
vertex $v_{1}$, both intersect all of the edges $v_{1}v_{l}$ exactly 
at $v_{1}$.  Both hyperplanes will miss all of the other edges 
$v_{i}v_{j}$, with $i,j \geq 2$, since for every point on all of these 
edges, 
$x_{1}=1$, $x_{2}=1$ and  $x_{l}\geq 0$, forcing the intersections to 
be empty.  But now we just have to count and see that the number of 
vertices is indeed $n+1$.  Thus $\bigtriangleup_{k}$ is a simplex.

The argument is similar for $\bigtriangleup_{ij}$.  Here we look at the 
hyperplanes
$$x_{1}+\ldots+ x_{i} +  x_{j} = 1$$
and 
$$x_{1}+\ldots+ x_{n-1} +  x_{j+1} = 1.$$
Both will intersect each of the edges $v_{0}v_{l}$, for $l\neq j$, in the 
same point, will intersect the edges  $v_{1}v_{l}$ exactly 
at $v_{1}$ and will miss the edges $v_{i}v_{j}$, with $i,j \geq 2$.
 They will intersect the edge $v_{0}v_{j}$ at distinct 
points.  Then $\bigtriangleup_{ij}$ has $n+1$ distinct 
vertices and is thus a simplex.

\noindent{\bf QED}

 Define the {\it n-triangle map}
 $T:\bigtriangleup \rightarrow \bigtriangleup $ by 
setting
$$T(\alpha_{1}, \ldots, \alpha_{n}) =
 (\frac{\alpha_{2}}{\alpha_{1}}, \ldots, 
 \frac{\alpha_{n-1}}{\alpha_{1}}, \frac{1 - \alpha_{1} \ldots 
 -\alpha_{n-1} - k 
\alpha_{n}}{\alpha_{1}}),$$
if  $(\alpha_{1}, \ldots, \alpha_{n}) \in \bigtriangleup_{k}$ and by
$$T(\alpha_{1}, \ldots, \alpha_{n}) =
 (\frac{\alpha_{2}}{\alpha_{1}}, \ldots, 
 \frac{\alpha_{j}}{\alpha_{1}}, \frac{1 - \alpha_{1} \ldots 
 -\alpha_{i}}{\alpha_{1}}, \frac{\alpha_{j+1}}{\alpha_{1}}, \ldots, 
 \frac{\alpha_{n}}{\alpha_{1}}),$$
if  $(\alpha_{1}, \ldots, \alpha_{n}) \in \bigtriangleup_{ij}$.

By direct calculation, we 
see that $T(\alpha_{1}, \ldots, \alpha_{n})\in\bigtriangleup$.  
Further,  each of the restriction maps
 $T:\bigtriangleup_{k} \rightarrow \bigtriangleup$ and 
  $T:\bigtriangleup_{ij} \rightarrow \bigtriangleup$  are 
one-to-one and onto, since the vertices of $\bigtriangleup_{k}$ and 
$\bigtriangleup_{ij}$
 map 
to the vertices of $\bigtriangleup$ and since lines map to lines.

We want to associate to each    $(\alpha_{1}, \ldots, \alpha_{n})$ 
in $\bigtriangleup$ an infinite sequence $(a_{0}, a_{1}, \ldots)$, 
where each $a_{k}$ is either a non-negative integer or a symbol 
$(ij)$, where $1\leq i \leq n-2$ and $i < j \leq n$.  If 
$T^{k}(\alpha_{1}, \ldots, \alpha_{n}) \in 
 \bigtriangleup_{l}$, set $a_{k} = l$ and
  if $T^{k}(\alpha_{1}, \ldots, \alpha_{n}) \in 
 \bigtriangleup_{ij}$, set $a_{k} = (ij)$.  Finally, if the $n$th 
 term for $T^{k}(\alpha_{1}, \ldots, \alpha_{n})$ is zero, stop.
 \begin{definition}
The n-triangle sequence  of  $(\alpha_{1}, \ldots, \alpha_{n})$ is the sequence
$(a_{1}, \ldots)$.
\end{definition}

We can also recursively define the triangle sequence as follows.  We 
want to define a sequence of $(n+1)$-tuples of nonincreasing positive 
numbers.  We will denote this sequence by $d_{1}(k),\ldots, 
d_{n+1}(k)$, for $k\geq 0$. Start with  
$$d_{1}(0) = 1, d_{2}(0) = 
\alpha_{1},\ldots, 
d_{n+1}(0)= \alpha_{n}.$$  Assume we have $d_{1}(k-1),\ldots, 
d_{n+1}(k-1)$. Define the symbol $a_{k}$ as follows.   If there is a 
nonnegative integer $l$ such that 
$$d_{1}(k-1) - d_{2}(k-1) - \ldots - d_{n}(k-1)
 - l d_{n+1}(k-1) \geq 0$$
 but 
$$d_{1}(k-1) - d_{2}(k-1) - \ldots - d_{n}(k-1)
 - (l+1) d_{n+1}(k-1) < 0,$$
 set 
$a_{k}= l$ and define
$$d_{1}(k) = d_{2}(k-1), \ldots, d_{n}(k) = d_{n+1}(k-1)$$
and
$$d_{n+1}(k) = d_{1}(k-1) - d_{2}(k-1) - \ldots - d_{n}(k-1)
 - l d_{n+1}(k-1).$$
 If no such integer exists, then there is a pair $(ij)$, with $1\leq 
 i \leq n-1$ and $i < j\leq n+1$ such that
 $$d_{j}(k-1) \geq d_{1}(k-1) - d_{2}(k-1) - \ldots - d_{i}(k-1)
> d_{j+1}(k-1).$$
In this case, define $a_{k}= (ij)$ and set
$$d_{1}(k) = d_{2}(k-1), \ldots, d_{j-1}(k) = d_{j}(k-1),$$
$$d_{j}(k) = d_{1}(k-1) - d_{2}(k-1) - \ldots - d_{i}(k-1)$$
and 
$$d_{j+1}(k) = d_{j+1}(k-1), \ldots, d_{n+1}(k) = d_{n+1}(k-1).$$

Now for the matrix version.  Let $(a_{1}, \ldots)$ be an 
n-triangle sequence for $(\alpha_{1}, \ldots, \alpha_{n})$.  If 
$a_{k}$ is a nonnegative integer, let 
$P_{k}$ be the $(n+1)\times (n+1)$

$$ \pmatrix{0& 0&  \dots & 0 &1\cr 1 & 0 & \dots & 0 &-1 \cr \; & \;\vdots 
 & \; & \; \cr 0 & \ldots & 1 & 0 & -1 \cr 
0 & \dots &  0  & 1 &-a_{k}\cr}.$$

If $a_{k}$ is the pair $(ij)$, let $P_{k}$ be the $(n+1)\times (n+1)$ 
matrix defined by:
$$(x_{1}, \ldots, x_{n+1})P_{k} =(x_{2},\ldots x_{j}, x_{1}-x_{2} 
-\ldots - x_{i}, x_{j+1},\ldots, x_{n+1}).$$

Then set $M_{k} = P_{1}\cdot P_{2} \cdots P_{k}$.  Note that $\det 
M_{k} = \pm 1$.  We have
$$(d_{1}(k), \ldots d_{n+1}(k)) = (1, \alpha_{1}, \ldots, \alpha_{n})M_{k}.$$ 
Set $M_{k} = (C_{1}(k),\ldots,C_{n+1}(k))$, where each $C_{m}(k)$ is a column 
vector of the matrix.  Then, if $a_{k}=l$, 
$$C_{m}(k) = C_{m+1}(k-1)$$
for $m\leq n$ and
$$C_{n+1}(k) = C_{1}(k-1) - C_{2}(k-1) -\ldots - C_{n}(k) - l 
C_{n+1}(k-1).$$
If $a_{k}=(ij)$, then, for $1 \leq m \leq j+1$,
$$C_{m-1}(k) = C_{m}(k-1),$$
$$C_{j+1}(k) = C_{1}(k-1) - C_{2}(k-1) -\ldots - C_{i+1}(k),$$
and for $j+1 \leq m \leq n+1,$
$$C_{m}(k) = C_{m}(k-1).$$

Each $C_{k}(m)$ can be viewed as an element of the integer lattice 
$Z^{n+1}$.  Then we have a method for producing elements of the 
integer lattice that approach the hyperplane
$$x_{0} + \alpha_{1} x_{1} + \ldots + \alpha_{n}x_{n} =0.$$

It is still unknown how to determine when an n-triangle sequence will 
uniquely determine an n-tuple
 $(\alpha_{1}, \ldots, \alpha_{n})\in\bigtriangleup$.  If we have 
 uniqueness, we strongly suspect that

 $$\alpha_{j} = \lim_{k\rightarrow \infty} 
 (-1)^{j}\frac{M_{k}(j1)}{M_{k}(11)},$$
 where $M_{k}(ij)$ denotes the determinant of the $n\times  n$ minor 
 of $M_{k}$ obtained by deleting the ith row and jth column.   The moral, but incorrect argument, is 
 the following.  First, since $\det M_{k}=\pm 1$, its 
 column vectors are linearly independent.  But also, the column 
 vectors are approaching the hyperplane
  whose normal vector is $(1, \alpha_{1}, \ldots, \alpha_{n})$.  Then 
  via standard arguments, the wedge product $C_{2}(k)\wedge 
  \ldots \wedge C_{n+1}(k)$  corresponds under duality to a vector 
  perpendicular to $C_{2}(k),
  \ldots, C_{n+1}(k)$ and by normalizing, will
   approach the vector $(1, \alpha_{1}, \ldots, \alpha_{n})$

  With reasonable conditions about uniqueness, we should have
  \begin{conjecture}
Let $0 \leq  \alpha_{n} \leq \ldots  \leq  \alpha_{1} < 1$ be 
an n-tuple  of numbers whose triangle 
sequence is $(k,k,k,\ldots)$.  Then $ \alpha_{j} = \alpha_{1}^{j}$ and 
$\alpha_{1}$
is a root of the algebraic equation
$$x^{n+1} + k x^{n} + x^{n-1} + \ldots  + x -1 = 0.$$
Further if $\alpha$ is the real root of this equation that is between
zero and one, then 
$(\alpha, \alpha^{2}, \ldots, \alpha^{n})$ has purely periodic simplex
 sequence $(k,k,k, 
\ldots)$. 
\end{conjecture}
A similar result holds if the triangle sequence is purely periodic 
of period length one of the form $(ij, ij, ij, \ldots)$.

We should also have
 \begin{conjecture}
 Let $(\alpha_{1}, \ldots, \alpha_{n})$ be an n-tuple real  numbers in 
 $\bigtriangleup$ whose triangle sequence is eventually periodic.  
 Then each $\alpha_{j}$  is algebraic of  degree at 
 most n.  Finally, as a vector space over ${\bf Q}$, the dimension of 
 ${\bf Q}[\alpha_{1}, \ldots, \alpha_{n}]$ is at most n.
 \end{conjecture}

Finally a comment about notation.  In the abstract it is claimed that 
that we will express each $n$-tuple will be associated to a sequence 
of nonnegative integers but in this section we look at sequences not 
just of nonnegative integers but also of terms of the form $(ij)$ with 
$ 1 \leq i \leq n-2$ and $i < j \leq n$.  But there are only a finite 
number ( $n(n-2) + \frac{(n-2)(n-1)}{2}$) of these extra symbols.  
We could, if desired, order these symbols by 
nonnegative numbers $0,1,\ldots, n(n-2) + \frac{(n-2)(n-1)}{2}$ and 
then shift the original nonnegative integers by this amount.  This 
will force the sequence to be one of nonnegative integers, but the 
notation is clearly worse than the one chosen.

\end{document}